\documentclass[11pt]{amsart}
\input{xy}
\xyoption{all}

\newcommand{\w}{\omega}
\newcommand{\IN}{\mathbb N}

\newcommand{\U}{\mathcal U}

\newtheorem{theorem}{Theorem}
\newtheorem{problem}{Problem}

\newtheorem{proposition}{Proposition}
\theoremstyle{definition}
\newtheorem{example}{Example}

\title{Characterizing meager paratopological groups}
\author{T.Banakh, I.Guran and A.Ravsky}
\address{Uniwersytet Humanistyczno-Przyrodniczy Jana Kochanowskiego w Kielcach, Poland}
\address{Ivan Franko National University of Lviv, Ukraine}
\email{T.O.Banakh@gmail.com, igor\_guran@yahoo.com, oravsky@mail.ru}


\begin{document}
\begin{abstract} We prove that a Hausdorff paratopological group $G$ is meager if and only if there are a nowhere dense subset $A\subset G$ and a countable set $C\subset G$ such that $CA=G=AC$.
\end{abstract}
\maketitle

Trying to find a counterpart of the Lindef\"of property in the category of topological groups I.Guran \cite{Gu} introduced the notion of an $\w$-bounded group which turned out to be very fruitful in topological algebra, see \cite{Tk}. We recall that a topological group $G$ is {\em $\omega$-bounded\/} if for each non-empty open subset $U\subset G$ there is a countable subset $C\subset G$ such that $CU=G=UC$.

A similar approach to the Baire category leads us to the notion of an
{\em shift-meager} ({\em shift-Baire}) group. This is a topological group that can(not) be written as the union of countably many translation copies of some fixed nowhere dense subset.

The notion of a shift-meager (shift-Baire) group can be defined in a more general context of semitopological groups, that is, groups $G$ endowed with a shift-invariant topology $\tau$. The latter is equivalent to saying that the group operation $\cdot:G\times G\to G$ is separately continuous. If this operation is jointly continuous, then $(G,\tau)$ is called a
{\em paratopological group}, see \cite{AT}.
\smallskip

A semitopological group $G$ is defined to be
\begin{itemize}
\item {\em left meager} (resp. {\em right meager}) if $G=CA$ (resp. $G=AC$) for some  nowhere dense subset $A\subset X$ and some countable subset $C\subset G$;
\item {\em shift-meager} if $G$ is both left and right meager;
\item {\em left Baire} (resp. {\em right Baire}) if for every open dense subset $U\subset X$ and every countable subset $C\subset G$ the intersection $\bigcap_{x\in C}xU$ (resp. $\bigcap_{x\in C}Ux$) is dense in $G$;
\item {\em shift-Baire} if $G$ is both left and right Baire.
\end{itemize}
For semitopological groups those notions relate as follows:

$$
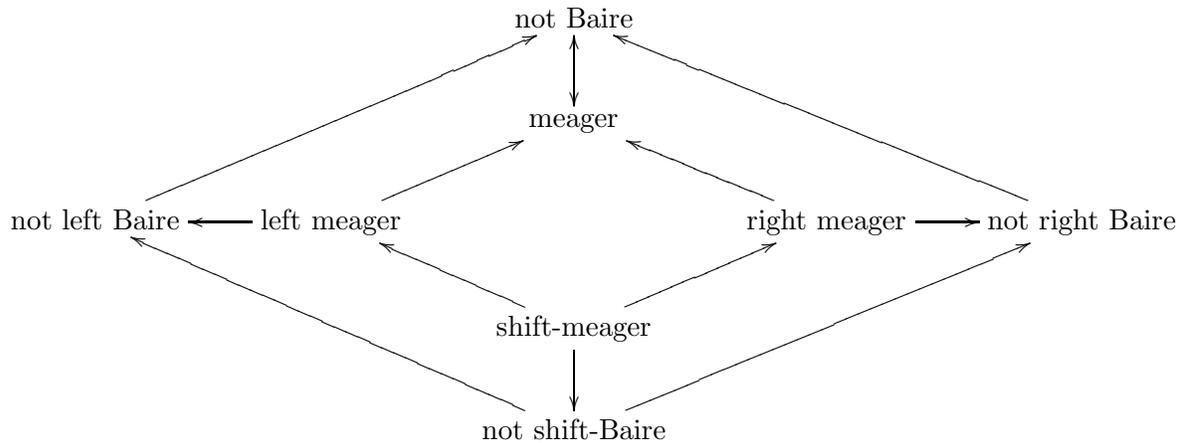
\xymatrix{
&&\mbox{not Baire}\\
&&\mbox{meager}\ar@{<->}[u]\\
\mbox{not left Baire}\ar[uurr]&\mbox{left meager}\ar[l]\ar[ur]& &\mbox{right meager}\ar[r]\ar[lu]&\mbox{not right Baire}\ar[lluu]\\
&&\mbox{shift-meager}\ar[lu]\ar[d]\ar[ur]\\
&&\mbox{not shift-Baire}\ar[lluu]\ar[rruu]
}
$$

The following theorem implies that for Hausdorff paratopological groups all the eight properties from this diagram are equivalent.

\begin{theorem}\label{main} A Hausdorff paratopological group $G$ is meager if and only if $G$ is shift-meager.
\end{theorem}

This theorem will be proved in Section~\ref{pf}. The proof is based on Theorem~\ref{t2} giving conditions under which a meager semitopological group is left (right) meager and Theorem~\ref{osc} describing some oscillator properties of 2-saturated Hausdorff paratopological groups.

\section{Shift-meager semitopological groups}

In this section we search for conditions under which a given meager semitopological group is left (right) meager.

Following \cite{BM1}, \cite{BM2} and \cite{DMM}, \cite{DP}, we define a subset $A\subset G$ of a group $G$ to be
\begin{itemize}
\item {\em left large} (resp. {\em right large}) if $G=FA$ (resp. $G=AF$) for some finite subset $F\subset G$;
\item {\em left P-small} (resp. {\em right P-small}\/) if there is an infinite subset $B\subset G$ such that the indexed family $\{bA\}_{b\in B}$ (resp. $\{Ab\}_{b\in B}$) is disjoint.
\end{itemize}

\begin{theorem}\label{t2} A meager semitopological group $G$ is left (right) meager provided one of the following conditions holds:
\begin{enumerate}
\item[(1)] $G$ contains a non-empty open left (right) P-small subset;
\item[(2)] $G$ contains a sequence $(U_n)_{n\in\w}$ of pairwise disjoint open left (right) large subsets;
\item[(3)] $G$ contains sequences of non-empty open sets  $(U_n)_{n\in\w}$ and points $(g_n)_{n\in\w}$ such that the sets $g_nU_nU_n^{-1}$ (resp. $U^{-1}_nU_ng_n$), $n\in\IN$, are pairwise disjoint.
\end{enumerate}
\end{theorem}

\begin{proof} $(1_l)$ Assume that $U\subset G$ is a non-empty open left P-small subset. We may assume that $U$ is a neighborhood of the neutral element $e$ of $G$.
It follows that there is a countable subset $B=\{b_n\}_{n\in\w}\subset G$ such that $b_nU\cap b_mU=\emptyset$ for any distinct numbers $n\ne m$. The countable set $B$ generates a countable subgroup $H$ of $G$. By an {\em $H$-cylinder} we shall understand an open  subset of the form $HVg$ where $g\in G$ and $V\subset U$ is a neighborhood of $e$.
Let $\U=\{HV_\alpha g_\alpha:\alpha\in A\}$ be a maximal disjoint family of $H$-cylinders in $G$ (such a family exists by the Zorn Lemma).

We claim that $\cup\U$ is dense in $G$. Assuming the converse, we could find a point $g\in G\setminus\cup\U$ and a neighborhood $V\subset U$ of $e$ such that $Vg\cap\cup\U$.
Taking into account that $H\cdot(\cup\U)=\cup\U$, we conclude that $HVg\cap\cup\U=\emptyset$ and hence $\U\cup\{HVg\}$ is a disjoint family of $H$-cylinders that enlarges the family $\U$, which contradicts the maximality of $\U$. Therefore $\cup\U$ is dense in $G$ and hence $G\setminus\cup\U$ is a closed nowhere  dense subset of $G$.

The space $G$, being meager, can be written as the union $G=\bigcup_{n\in\w}M_n$ of a  sequence $(M_n)_{n\in\w}$ of nowhere dense subsets of $G$. It is easy to see that the set
$$M=(G\setminus\cup\U)\cup\bigcup_{\alpha\in A}\bigcup_{n\in\w}b_n(M_n\cap V_\alpha g_\alpha)$$
is nowhere dense in $G$ and $G=HM$, witnessing that $G$ is left meager.
\smallskip

$(2_l)$ Assume that $G$ contains a sequence $(U_n)_{n\in\w}$ of pairwise disjoint open left large subsets. For every $n\in\w$ find a finite subset $F_n\subset G$ with $G=F_n\cdot U_n$.
Write $G=\bigcup_{n\in\w}M_n$ as countable union of nowhere dense subsets and observe that for every $n\in\w$ the subset $\bigcup_{x\in F_n}x^{-1}(M_n\cap x U_n)$ of $U_n$ is nowhere dense. Since the family $\{U_n\}_{n\in\w}$ is disjoint, the set
$$M=\bigcup_{n\in\w}\bigcup_{x\in F_n}x^{-1}(M_n\cap xU_n)$$is nowhere dense in $G$.
Since $(\bigcup_{n\in\w}F_n)\cdot M=G$, the semitopological group $G$ is left meager.
\smallskip

$(3_l)$ Assume that $(U_n)_{n\in\w}$ is a sequence of non-empty open subsets of $G$ and $(g_n)_{n\in\w}$ is a sequence of points of $G$ such that the sets  $g_nU_n U_n^{-1}$, $n\in\w$, are pairwise disjoint. Using the Zorn Lemma, for every $n\in\w$ we can choose a maximal subset $F_n\subset G$ such that the indexed family $\{x U_n\}_{x\in F_n}$ is disjoint. If for some $n\in\w$ the set $F_n$ is infinite, then the set $U_n$ is left P-small and consequently, the group $G$ is left meager by the first item. So, assume that each set $F_n$, $n\in\w$, is finite. The maximality of $F_n$ implies that for every $x\in G$ there is $y\in F_n$ such that $xU_n\cap yU_n\ne\emptyset$. Then $x\in yU_nU_n^{-1}$ and hence $G=F_nU_nU_n^{-1}$, which means that the open  set
$U_nU_n^{-1}$ is left large. Since the family $\{g_nU_n U_n^{-1}\}_{n\in\w}$ is disjoint, it is legal to apply the second item to conclude that the group $G$ is left meager.
\smallskip

$(1_r)-(3_r)$. The right versions of the items (1)--(3) can be proved by analogy.
\end{proof}

\section{Oscillation properties of paratopological groups}

In this section we establish some oscillation properties of 2-saturated paratopological groups. First, we recall the definition of oscillator topologies on a given paratopological group $(G,\tau)$, see \cite{BR} for more details.

Given a subset $U\subset G$, by induction define subsets
$(\pm U)^n$ and $(\mp U)^n$, $n\in\omega$, of $G$ letting $(\pm U)^0=(\mp
U)^0=\{e\}$ and $(\pm U)^{n+1}=U(\mp U)^n$, $(\mp
U)^{n+1}=U^{-1}(\pm U)^n$ for $n\ge 0$. Thus $(\pm
U)^n=\underset{n}{\underbrace{UU^{-1}U\cdots U^{(-1)^{n-1}}}}$ \
and $(\mp U)^n=\underset{n}{\underbrace{U^{-1}UU^{-1}\cdots
U^{(-1)^{n}}}}$. Note that $((\pm U)^n)^{-1}=(\pm U)^n$ if $n$ is
even and $((\pm U)^n)^{-1}=(\mp U)^n$ if $n$ is odd.

By an {\em $n$-oscillator\/} (resp. {\em a mirror
$n$-oscillator\/}) on a topological group $(G,\tau)$ we understand
a set of the form $(\pm U)^n$ (resp. $(\mp U)^n$~) for some
neighborhood $U$ of the unit of $G$. Observe that each
$n$-oscillator in a paratopological group $(G,\tau)$ is a mirror
$n$-oscillator in the mirror paratopological group $(G,\tau^{-1})$
and vice versa: each mirror $n$-oscillator in $(G,\tau)$ is an
$n$-oscillator in $(G,\tau^{-1})$.

By the {\em $n$-oscillator topology\/} on a paratopological
group $(G,\tau)$ we understand the topology $\tau_n$ consisting of
sets $U\subset G$ such that for each $x\in U$ there is an
$n$-oscillator $(\pm V)^n$ with $x\cdot(\pm V)^n\subset U$.

Let us recall \cite{Gu2}, \cite[p.342]{AT} that a paratopological group $(G,\tau)$ is {\em saturated} if each non-empty open set $U\subset G$ has non-empty interior in the mirror topology $\tau^{-1}=\{U^{-1}:U\in\tau\}$. This notion can be generalized as follows.

Define a paratopological group $(G,\tau)$ to be {\em $n$-saturated} if each non-empty open set $U\in\tau_n$ has non-empty interior in the topology $(\tau^{-1})_n$.

\begin{proposition}\label{p1} A paratopological group $(G,\tau)$ is 2-saturated if no non-empty open subset $U\subset G$ is P-small.
\end{proposition}

\begin{proof} To prove that $G$ is 2-saturated, take any non-empty open set $U_2\in\tau_2$ and find a point $x\in U_2$ and a neighborhood $U\in\tau$ of $e$ such that $xU^2U^{-2}\subset U_2$. By the Zorn Lemma, there is a maximal subset $B\subset G$ such that $bU\cap b'U=\emptyset$ for all distinct points $b,b'\in B$. By our hypothesis, $U$ is not P-small, which implies that the set $B$ is finite.
The maximality of $B$ implies that for each $x\in G$ the shift $xU$ meets some shift $bU$, $b\in B$. Consequently, $x\in bUU^{-1}$ and $G=\bigcup_{b\in B}bUU^{-1}$.
It follows that the closure $\overline{UU}^{-1}$ of $UU^{-1}$ in the topology $(\tau^{-1})_2$  has non-empty interior. We claim that  $\overline{UU}^{-1}\subset U^2U^{-2}$. Indeed, given any point $z\in \overline{UU}^{-1}$, we conclude that the neighborhood $U^{-1}zU$ of $z$ in the topology $(\tau^{-1})_2$ meets $UU^{-1}$ and hence $z\in U^2U^{-2}$. Now we see that the set $$U_2\supset xU^2U^{-2}\supset x\overline{UU^{-1}}$$ has non-empty interior in the topology $(\tau^{-1})_2$, witnessing that the group $(G,\tau)$ is 2-saturated.
\end{proof}

By Proposition 2 of \cite{BR}, for each saturated paratopological group $(G,\tau)$ the semitopological group $(G,\tau_2)$ is a topological group. This results generalizes to $n$-saturated groups.

\begin{theorem}\label{osc} If $(G,\tau)$ is an $n$-saturated paratopological group for some $n\in\IN$, then $(G,\tau_{2n})$ is a topological group.
\end{theorem}

\begin{proof} According to Theorem 1 of \cite{BR}, $(G,\tau_{2n})$ is a topological group if and only if for every neighborhood $U\in\tau$ of the neutral element $e\in G$ there is a neighborhood $V\in\tau$ of $e$ such that $(\mp V)^{2n}\subset (\pm U)^{2n}$.

Since the paratopological group $(G,\tau)$ is $n$-saturated, the set $(\pm U)^n\in\tau_2$ contains an interior point $x$ in the mirror topology $(\tau^{-1})_n$. Consequently, there is a neighborhood $V\in\tau$ of $e$ such that $(\mp V)^n x\subset (\pm U)^n$.

Now we consider separately the cases of odd and even $n$.

1. If $n$ is odd, then applying the operation of the inversion to $(\mp V)^n x\subset (\pm U)^n$, we get $x^{-1}(\pm V)^n\subset (\mp U)^n$ and then
$$(\mp V)^{2n}=(\mp V)^n (\pm V)^n=(\mp V)^nxx^{-1}(\pm V)^n\subset (\pm U)^n(\mp U)^n=(\pm U)^{2n}.$$

2. If $n$ is even, then $(\mp V)^n x\subset (\pm U)^n$ implies $x^{-1}(\mp V)^n\subset (\pm U)^n$ and
$$(\mp V)^{2n}=(\mp V)^n (\mp V)^n=(\mp V)^nxx^{-1}(\mp V)^n\subset (\pm U)^n(\pm U)^n=(\pm U)^{2n}.$$
\end{proof}

According to \cite{BR}, for each 1-saturated Hausdorff paratopological groups $(G,\tau)$ the group $(G,\tau_2)$ is a Hausdorff topological group. For 2-saturated group we have a bit weaker result.

\begin{theorem}\label{t4} For any non-discrete Hausdorff 2-saturated paratopological group $(G,\tau)$  the maximal antidiscrete subgroup $\overline{\{e\}}=\bigcap_{e\in U\in\tau}(\pm U)^4$ of the topological group $(G,\tau_4)$ is nowhere dense in the topology $\tau_2$.
\end{theorem}

\begin{proof} To show that  $\overline{\{e\}}$ is nowhere dense in the topology $\tau_2$, fix any non-empty open set $U_2\in\tau_2$. Since $G$ is not discrete, so is the topology $\tau_2\subset\tau$. Consequently, we can find a point $x\in U_2\setminus\{e\}$. Since $G$ is a Hausdorff paratopological group, there is a neighborhood $U\in\tau$ of $e$ such that $e\notin xUU^{-1}\subset U_2$. The continuity of the group operation yields a neighborhood $V\in\tau$ of $e$ such that $V^2\subset U$ and $V^2x\subset xU$. Then $V^2xV^{-2}\subset xUU^{-1}\not\ni e$ yields  $V^{-1}V\cap VxV^{-1}=\emptyset$.
Using the shift-invariantness of the topology $\tau$, find a neighborhood $W\in\tau$ of $e$ such that $W\subset V$ and $xW\subset Vx$.

Since the group $G$ is 2-saturated, the open set $xWW^{-1}\in\tau_2$ has non-empty interior in the topology $(\tau^{-1})_2$. Consequently, there is a point $y\in xWW^{-1}$ and a neighborhood $O\in\tau$ of $e$ such that $O\subset W$ and $O^{-1}yO\subset xWW^{-1}$. Observe that
$$O^{-1}O\cap O^{-1}yO\subset V^{-1}V\cap xWW^{-1}\subset V^{-1}V\cap VxV^{-1}=\emptyset$$ and consequently,
$U_2\ni y\notin OO^{-1}OO^{-1}\supset\overline{\{e\}}$.
\end{proof}

\begin{problem} Can the topology $\tau_{2n}$ be antidiscrete for some Hausdorff $n$-saturated paratopological group?
\end{problem}

\begin{problem} Assume that a paratopological group $(G,\tau)$ is 2-saturated. Is its mirror paratopological group $(G,\tau^{-1})$ 2-saturated?
\end{problem}

\section{Proof of Theorem~\ref{main}}\label{pf}

We need to check that each meager Hausdorff paratopological group $G$ is left and right meager.

If the paratopological group $G$ contains a non-empty open left P-small subset, then $G$ is left meager by Theorem~\ref{t2}(1). So assume that no non-empty open subset of $G$ is left P-small. In this case Proposition~\ref{p1} implies that the paratopological group $G$ is 2-saturated while Theorem~\ref{t4} ensures that the topological group $(G,\tau_4)$ contains a countable disjoint family  $\{W_n\}_{n\in\w}$ of non-empty open sets. By the definition of the 4th oscillator topology $\tau_4$ each set $W_n$ contains a subset of the form $x_nU_nU^{-1}U_nU_n^{-1}$ where $x_n\in G$ and $U_n$ is a neighborhood of the neutral element in the paratopological group $G$. Since the sets $x_nU_nU_n^{-1}\subset W_n$, $n\in\w$, are pairwise disjoint, we can apply Theorem~\ref{t2}(3) to conclude that the paratopological group $G$ is left meager.

By analogy we can prove that $G$ is right meager.

\section{Discussion and Open Problems}

The following example shows that without any restrictions, a meager semi-topological group needs not be shift-meager.

\begin{example} Let $G$ be an uncountable group whose cardinality $|G|$ has countable cofinality. Endow the group $G$ with the shift-invariant topology generated by the base $\{G\setminus A:|A|<|G|\}$. It is easy to see that a subset $A\subset G$ is nowhere dense if and only if it is not dense if and only if $|A|<|G|$. This observation implies that $G$ is meager (because $|G|$ has countable cofinality). On the other hand, the semi-topological group $G$ is not shift-meager because for every nowhere dense subset $A\subset G$ and every countable subset $C\subset G$ we get $|A|<|G|$ and hence $G\ne CA$ because $|CA|\le\max\{\aleph_0,|A|\}<|G|$.
\end{example}

\begin{problem} Is each meager paratopological group $G$ shift-meager?
\end{problem}

\begin{problem} Is each meager Hausdorff semitopological group shift-meager?
\end{problem}

\begin{problem}
Is each left meager semitopological group right meager?
\end{problem}

Also we do not know is the following semigroup version of Theorem~\ref{main} holds.

\begin{problem} Let $S$ be an open meager subsemigroup of a Hausdorff paratopological group $G$. Is $S\subset CA$ for some nowhere dense subset $A\subset S$ and a countable subset $C\subset G$?
\end{problem}

\end{document}